\documentclass[11pt]{amsart}
\usepackage{amsmath,amssymb,latexsym}
\usepackage[T1]{fontenc}
\usepackage{enumerate}

\def\M{\mathfrak M}

\def\Aut{\mathrm{Aut}}

\def\tp{\mathrm{tp}}
\def\lstp{\mathrm{Lstp}}
\def\stp{\mathrm{stp}}

\def\SU{\mathrm{SU}}
\def\cb{\mathrm{Cb}}
\def\acl{\mathrm{acl}}
\def\dcl{\mathrm{dcl}}
\def\bdd{\mathrm{bdd}}

\def\Ind#1#2{#1\setbox0=\hbox{$#1x$}\kern\wd0\hbox to 0pt{\hss$#1\mid$\hss}
\lower.9\ht0\hbox to 0pt{\hss$#1\smile$\hss}\kern\wd0}
\def\ind{\mathop{\mathpalette\Ind{}}}
\def\Notind#1#2{#1\setbox0=\hbox{$#1x$}\kern\wd0\hbox to 0pt{\mathchardef
\nn="3236\hss$#1\nn$\kern1.4\wd0\hss}\hbox to 0pt{\hss$#1\mid$\hss}\lower.9\ht0
\hbox to 0pt{\hss$#1\smile$\hss}\kern\wd0}
\def\nind{\mathop{\mathpalette\Notind{}}}

\theoremstyle{plain}
\newtheorem{theorem}{Theorem}[section]
\newtheorem{prop}[theorem]{Proposition}
\newtheorem{fact}[theorem]{Fact}
\newtheorem{lemma}[theorem]{Lemma}
\newtheorem{cor}[theorem]{Corollary}

\theoremstyle{definition}
\newtheorem{defn}[theorem]{Definition}
\newtheorem{remark}[theorem]{Remark}

\newtheorem{expl}[theorem]{Example}
\newtheorem*{conj}{Conjecture}

\def\pf{\par\noindent{\em Proof. }}

\title[Generalized amalgamation and homogeneity]{Generalized amalgamation and homogeneity}
\author{Daniel Palac\'in}
\thanks{This work was partially supported by the project MTM2014-59178-P and the project SFB 878.
I am indebted with Vera Koponen for pointing out a mistake in an earlier version of the paper and also I would like to thank her for several discussions.}
\address{Mathematisches Institut, Universit\"at M\" unster, Einsteinstrasse 62,
48149 M\"unster, Germany}
\email{daniel.palacin@uni-muenster.de}

\begin{document}

\begin{abstract}
In this paper we shall prove that any $2$-transitive finitely homogeneous structure with a supersimple theory satisfying a generalized amalgamation property is a random structure. In particular, this adapts a result of Koponen for binary homogeneous structures to arbitrary ones without binary relations. Furthermore, we point out a relation between generalized amalgamation, triviality and quantifier elimination in simple theories.
\end{abstract}

\maketitle

\section{Introduction}

A permutation group acting on a set $X$ is said to be {\em oligomorphic} if its action has only finitely many orbits on $X^n$ for each natural number $n$. Such groups appear naturally in model theory as groups of automorphisms of
countable $\omega$-categorical structures.  In fact, after a theorem of Ryll-Nardzewski, any countable structure $\M$ (in a countable language) is $\omega$-categorical if and only if its group of automorphisms $\Aut(\M)$ is oligomorphic. Within the family of $\omega$-categorical structures lies the large class of {\em finitely homogeneous} relational ones. These are those countable structures in a finite relational language such that any isomorphism between finite substructures can be extended to an automorphism of the whole structure. Alternatively, finitely homogeneous structures are precisely those structures in a finite relational language that can be built up from its finite substructures by the Fra\"iss\'e amalgamation method. 
Classical examples are countable abelian groups of finite exponent, the countable dense linear order without endpoints and the random graph.

From the point of view of classification theory, $\omega$-categorical structures form an important class. For instance, the random graph, or more generally a random structure in the sense of Definition \ref{DefRandom}, is an archetypical example of a structure with a first-order simple, even supersimple, theory. In these structures, there is always a well-behaved notion of independence among subsets, called {\em forking independence}, which satisfies certain amalgamation property known as the Independence Theorem or also $3$-complete amalgamation, see Section 2. In the case of the random graph this independence is given by equality, i.e. a set $A$ is independent from $B$ over $C$ if and only if $A\cap B$ is contained in $C$. Moreover, it satisfies the $n$-complete amalgamation property (see Definition \ref{DefAmalg}) for any natural number $n$ since any finite graph is embeddable in the random graph. However, not all simple structures do. For instance, the tetrahedron-free ternary random graph is simple but does not satisfy the $4$-complete amalgamation property. Thus, there is an intuitive connection between amalgamation properties and the existence of forbidden substructures which we shall formalize in Section $4$ in the context of finitely homogeneous structures.

We focus our attention on structures where $\Aut(\M)$ acts $2$-transitively on $\M$, which extends the notion of {\em primitive} structure, i.e. those structures $\M$ where there is no equivalence relation in $\M\times\M$ which is invariant under $\Aut(\M)$. The random graph is an example of a primitive supersimple structure, and in fact Koponen \cite{Koponen2} conjectured that it is essentially the canonical one. More precisely:

\begin{conj}
Any primitive binary finitely homogeneous supersimple structure is a random structure.
\end{conj}

In \cite{Koponen3}, Koponen solved the conjecture under the assumption of {\em one-basedness} (see Section 2) by showing that such structures have $\SU$-rank one and hence are random by a result of Aranda L\'opez \cite{Andres}.  In fact, an easy argument contained in \cite{Koponen1} yields that any binary finitely homogeneous supersimple (even simple) structure is one-based and so the full conjecture follows. In this paper we analyze the corresponding conjecture for arbitrary languages and show that $2$-transitive finitely homogeneous supersimple structures in a relational language with relations of arity at most $n$ are random whenever the $(n+1)$-complete amalgamation property holds (Theorem \ref{ThmRandom}). In particular, this applies to finitely homogeneous structures without binary relations. Additionally, if the theory eliminates imaginaries geometrically then the same is true for primitive structures (Remark \ref{RemPrimitive}). Finally, observe that as any simple structure satisfies $3$-complete amalgamation, in the binary case no additional assumptions are necessary. However, for relational languages of arity $n$, it is essential to assume $(n+1)$-complete amalgamation, as exhibited by the tetrahedron-free ternary random graph.

Our proof boils down to the study of forking independence in simple structures whose theory admits quantifier elimination in a relational language, where there is a bound on the arity of all relations. This is treated in full generality in Section $3$, where distinct degrees of triviality are pointed out under the assumption of quantifier elimination and complete amalgamation. Additionally, some examples, such as right-angled buildings with infinite residues \cite{BMPZ}, are provided to exemplify that the remarks in that section are interesting by themselves. In particular, it follows that in order to construct a non-trivial superstable (even supersimple) structure satisfying certain geometrical properties, which is a major problem in model theory, one should consider richer languages.

Throughout the text, we assume that the reader is familiar with the basics of model theory and simple theories. Nevertheless, in Section $2$, we recall  some of the basic notions of simplicity theory and also give the definition of the $n$-complete amalgamation property. For further details we refer the interested reader to \cite{cas-book,Kim,wagner}.

\section{Preliminaries on simplicity}

We shall be working inside a large saturated and homogeneous model $\M$ of a complete first-order theory. Thus, tuples and sets consist of elements from this large model.

We introduce simplicity in terms of dividing. A partial type $\pi(x)$ {\em divides} over $A$ if it implies a formula $\varphi(x,a)$ for which there is an $A$-indiscernible sequence $(a_i)_{i<\omega}$ in $\tp(a/A)$ such that $\{\varphi(x,a_i)\}_{i<\omega}$ is inconsistent. The theory is {\em simple} if any complete type $\tp(a/B)$ does not divide over a subset $A$ of $B$ of size $|A|\le |T|$. Additionally, it is said to be {\em supersimple} if $A$ can be taken to be finite. As a consequence, one can observe that the imaginary expansion of a (super)simple theory is again (super)simple. Thus, there is no harm in assuming that our large model is a model of the imaginary expansion of the theory.

In simple theories, dividing agrees with the notion of {\em forking}: A partial type forks over a set $A$ if it implies a finite disjunction of formulas, each of which divides over $A$.
Both notions give rise to a notion of independence among subsets of our model. Namely, we say that a set $A$ is independent from $B$ over $C$ if and only if the type $\tp(a/BC)$ of any finite tuple $a$ of elements from $A$ over $B\cup C$ does not fork over $C$. We write $A\ind_C B$ for this. Forking independence in simple theories is a well-behaved notion of independence, see \cite[Chapter 12]{cas-book} for an abstract approach to independence relations. Here we summarize some of its main properties:
\begin{enumerate}
\item Invariance under $\Aut(\M)$.
 \item Finite character. $A\ind_C B$ if and only if $A'\ind_C B'$ for any finite $A'\subseteq A$ and $B'\subseteq B$.
 \item Symmetry. If $A\ind_C B$, then $B\ind_C A$.
 \item Transitivity. $A\ind_C B$ and $A\ind_{BC} D$ if and only if $A\ind_C BD$.
 \item Extension. If $A\ind_C B$, then for any $D$ there is some $A'\equiv_{BC} A$ with $A'\ind_C BD$.
 \item Local character. For any finite $A$ and any $B$, there exists a subset $C\subseteq B$ with $|C|\le |T|$ such that $A\ind_{C} B$.
\item Algebraicity: If $A\ind_C A$, then $A\subseteq\acl^{\rm eq}(A)$.
\end{enumerate}
Given two types $p$ and $q$ we say that $p$ is a {\em non-forking extension} of $q$ if $p$ extends $q$ and it does not fork over the parameters of $q$. Using this notion, we can define {\em stable} theories as those simple theories where any type over a model has a unique non-forking extension over a larger set of parameters. One of the main features of forking independence in arbitrary simple theories is the Independence Theorem for types over models which can be seen as a weakening of the uniqueness of non-forking extensions in stable theories: If $p_0(x)\in S(M)$ is a type over a model $M$ and $p(x)\in S(A)$ and $q(x)\in S(B)$ are non-forking extensions of $p_0$ such that $A$ is independent from $B$ over $M$, then $p(x)\cup q(x)$ does not fork over $M$. In other words, it is possible to amalgamate non-forking extensions of a common type based over a model as soon as the parameters are independent over this model. 
This motivates the definition of a {\em Lascar strong type}: A type over which the Independence Theorem holds. We denote the Lascar strong type of $a$ over $A$ by $\lstp(a/A)$. In $\omega$-categorical simple theories and in supersimple theories any type over an algebraically closed set of imaginaries is Lascar strong. We can associate to each Lascar strong type $p$ a minimal definably closed element $e$ such that $p_{\upharpoonright e}$ is still Lascar strong and $p$ does not fork over $e$. Such an element $e$ is called the canonical base of $p$ and we write $\cb(a/A)$ to denote the canonical base of $\lstp(a/A)$. However, a canonical base is a {\em hyperimaginary}, i.e. the equivalence class of a tuple (possibly infinite) of an $\emptyset$-type-definable equivalence relation. Thus, the general theory of simplicity must be developed in the context of hyperimaginaries, and hence it becomes more technical. As hyperimaginaries do not appear explicitly in our paper, we omit to treat them in this introductory  exposition, but we refer the interested reader to \cite[Chapters 15-17]{cas-book}.

Supersimple theories allow an ordinal-valued rank which is compatible with forking independence. This is the $\SU$-rank, which is the least function assigning to each complete type an ordinal such that $\SU(p)\ge \alpha+1$ if and only if there is a forking extension $q$ of $p$ with $\SU(q)\ge \alpha$. That is, the $\SU$-rank corresponds to the fundamental rank of forking among complete types. We say that a simple theory has finite $\SU$-rank or $\SU$-rank one when all its $1$-types do.

Two types are said to be {\em orthogonal} if any two realizations (of non-forking extensions) are independent, and a type $p\in S(A)$ is {\em regular} if it is orthogonal to its forking extensions, i.e. if for any $B\supseteq A$ and elements $a,b$ realizing $p$ with $a\ind_A B$ and $b\nind_A B$ we have that $a\ind_B b$. In a supersimple theory, every non-algebraic type is non-orthogonal to a regular type. Similarly, a non-algebraic type of finite $\SU$-rank is non-orthogonal to a type of $\SU$-rank one. Consequently, as it is exhibited along the paper, types of $\SU$-rank one in a theory of finite $\SU$-rank coordinate the whole structure, as do regular types in the supersimple case. Finally, recall that in a simple theory a partial $\pi(x)$ with parameters over a set $A$ is {\em one-based} if for any tuple $\bar a$ of realizations of $\pi$ and any sets $C\supseteq B\supseteq A$ we have that $\cb(\bar a/B)$ is contained in $\cb(\bar a/C)$. In a supersimple theory (or more generally, in a simple theory where hyperimaginaries are eliminable) we can reformulate one-basedness of $\pi$ as follows: for any tuple $\bar a$ of realizations of $\pi$ and any set $B\supseteq A$ we have that $\bar a$ is independent from $B$ over $\acl^{\rm eq}(A,\bar a)\cap\acl^{\rm eq}(B)$.

We introduce a definition of generalized amalgamation, which requires the use of hyperimaginary elements. However, in the $\omega$-categorical or the supersimple cases one can simply work with $\acl^{\rm eq}$ instead of $\bdd$.

\begin{defn}
Let $W$ be a collection of subsets of $\mathcal P(n)$ closed under subsets. We say that a family $\{p_s(x_s)\}_{s\in W}$ of types over set $A$ is an {\em independent system of types} if the following three conditions hold:
\begin{enumerate}
\item for any $s\subseteq t$ in $W$ we have that $x_s\subseteq x_t$ and $p_s(x_s)\subseteq p_t(x_t)$.
\end{enumerate}
If $a_s$ realizes $p_s$ then for $s\in W$:
\begin{enumerate}
\item[(2)] the tuple $(a_{\{i\}})_{i\in s}$ is independent over $A$, and
\item[(3)] we have that $a_s=\bdd(A,(a_{\{i\}})_{i\in s})$.
\end{enumerate}
\end{defn}

For a given natural number $n$, let $\mathcal P^-(n)=\mathcal P(n)\setminus\{n\}$, i.e. $\mathcal P^{-}(n)$ is the collection of all subsets of $\{0,\ldots,n-1\}$ except $\{0,\ldots,n-1\}$.

\begin{defn}\label{DefAmalg}
A theory has {\em $n$-complete amalgamation} if for any independent system of Lascar strong types $\{p_s(x_s)\}_{s\in\mathcal P^-(n)}$ with parameters over some set $A$ there is a type $p_n(x_n)$ over $A$ such that $\{p_s(x_s)\}_{s\in\mathcal P(n)}$ is also an  independent system of types over $A$.
Furthermore, we say that a theory satisfies the $n$-complete amalgamation over models if the above holds for independent systems of types over models.
\end{defn}

It is easy to see that the Independence Theorem corresponds to $3$-complete amalgamation and hence, any simple theory has $3$-complete amalgamation. Additionally, $n$-complete amalgamation implies $m$-complete amalgamation for $m\le n$.  Any stable theory has $n$-complete amalgamation over models but not necessarily over arbitrary boundedly closed sets. Moreover, it is easy to find simple theories without $n$-complete amalgamation for any natural number $n$ as it is exhibited in the following example.

\begin{expl} Let $2\le n<k$. We say that an $n$-graph $(V,R)$ is $k$-free if $V$ is an infinite set and $R$ is an $n$-ary relation on $V$ such that $R(a_1,\ldots,a_n)$ with $a_i\in V$ implies that $a_1,\ldots,a_n$ are distinct and there is no $k$-clique, i.e. there is no subset of $V$ with every subset of size $k$ satisfying $R(x_1,\ldots,x_n)$.

There is a unique countable $k$-free $n$-graph which is a primitive finitely homogeneous structure. It can be easily obtained by free amalgamation. Moreover, for any $n$ and $k\geq 3$ its theory is one-based simple unstable of $\SU$-rank one. It is easy to see that the theory of a $k$-free $n$-graph does not satisfies the $k$-complete amalgamation property. Thus, the triangle-free random graph (i.e. $n=2$ and $k=3$) is not simple. See also \cite{Hru}.
\end{expl}



To conclude this preliminary section, we show that many finitely homogeneous structures with a simple theory, such as the countable $k$-free $n$-graph with $k\ge 3$, have $\SU$-rank one and are one-based. In fact, we prove a more general statement which we believe might be folklore but we cannot it in the literature. Before, recall that a theory has {\em geometric elimination of imaginaries} if every imaginary is interalgebraic with a real tuple. Moreover, a theory has {\em weak elimination of imaginaries} if for any imaginary $e$ there is a real tuple $a$ such that $a\subseteq\acl(e)$ and $e\in\dcl^{\rm eq}(a)$.

\begin{prop}
Assume that a countable $\omega$-categorical structure has a disintegrated algebraic closure (i.e. the algebraic closure of a set equals to the union of the algebraic closure of its singletons). If its theory is simple and has geometric elimination of imaginaries, then it is supersimple and one-based.
\end{prop}
\pf We shall work inside a large saturated and homogeneous model of the theory. It is suffices to show that a set $A$ is independent from $B$ over the set $\acl(A)\cap\acl(B)$. If not, we can find some finite tuples $a$ and $b$ such that $\tp(a/b)$ forks over $\acl(a)\cap\acl(b)$. As the theory is $\omega$-categorical, the canonical base $\cb(a/b)$ is a single imaginary. Thus, by geometric elimination of imaginaries, there is no harm in assuming that $\cb(a/b)$ is a finite real tuple, and so there is some $c\in\cb(a/b)$ such that $c\not\in\acl(a)\cap\acl(b)$. However, this element $c$ is algebraic over some (any) Morley sequence $(a_i)_{i<\omega}$ in $\lstp(a/b)$ and so $c\in\acl(a_i)$ since the algebraic closure is disintegrated by assumption. Thus, as $c\in\acl(b)$ we obtain that $c\in\acl(a)$ by indiscernibility, a contradiction. \qed

Observe that in the above result $\omega$-categoricity is not necessary as it can be replaced by the condition that canonical bases (over finite tuples) be imaginaries. Now, as a consequence we obtain the following result due to Conant \cite[Corollary 7.14]{Con}:

\begin{cor}
A finitely homogeneous structure with a simple theory obtained by Fra\"iss\'e construction with free amalgamation is supersimple of $\SU$-rank one and one-based.
\end{cor}
\pf Any finitely homogeneous structure obtained via a Fra\"iss\'e construction with free amalgamation has trivial algebraic closure and moreover, it weakly eliminates imaginaries \cite[Lemma 2.7]{MT}. Thus, the statement follows by the previous result.\qed

\section{Triviality}

In this section we establish connections between several notions of triviality introduced in \cite{Poi}, complete amalgamation and elimination of quantifiers.

\begin{defn}
A theory is {\em $k$-trivial} if for every set $A$ and any $(k+1)$-independent tuples $a_0,a_1,\ldots,a_{k+1}$ over $A$ are independent over $A$. When $k=1$ we simply say that the theory is trivial.
\end{defn}

Goode \cite[Lemma 1]{Poi} proved that if a stable theory is trivial, then so is its imaginary expansion. The same proof yields:

\begin{fact}\label{FactTrivImag}
If a simple theory is $k$-trivial among real elements, then so is among hyperimaginaries.
\end{fact}

\begin{prop}\label{PropTriv}
A simple theory with elimination of quantifiers in a relational language of arity at most $k+1$ is $k$-trivial whenever it has $(k+2)$-complete amalgamation over models.
\end{prop}
\pf By Fact \ref{FactTrivImag}, it is enough to show that the theory is $k$-trivial among real tuples.
Let $a_0,a_1,\ldots,a_{k+1}$ be finite tuples from the home sort and suppose that they are $(k+1)$-independent over a set $A$. Consider a model $M$ containing $A$ independent from $a_0,\ldots,a_{k+1}$ over $A$. Thus $a_0,\ldots,a_{k+1}$ are $(k+1)$-independent over the model $M$.

Now, for a subset $s$ from $\mathcal P^-(k+2)$ set $b_s=\bdd(M,(a_j)_{j\in s})$ and $p_s(x_s)=\tp(b_s/M)$. It is clear that  $x_s$ is contained in $x_t$ and  $p_s\subseteq p_t$  when $s\subseteq t$ for any $s,t\in \mathcal P^{-}(k+2)$.  Moreover, by construction
$$
b_s=\bdd(M,(b_{\{j\}})_{j\in s})
$$
and additionally, since $a_0,\ldots,a_{k+1}$ are $(k+1)$-independent over $M$, then so are $b_{\{0\}},\ldots,b_{\{k+1\}}$. Thus, the family $\{p_s(x_s)\}_{s\in\mathcal P^{-}(k+2)}$ is an independent system of types over $M$ and hence, by $(k+2)$-complete amalgamation over models we get a complete type $p_{k+2}(x_{k+2})$ over $M$ such that $\{p_s(x_s)\}_{s\in\mathcal P(k+2)}$ is an independent system of types as well.

Let $c$ be realization $p_{k+2}(x_{k+2})$ and for $s\in \mathcal P^{-}(k+2)$, set $c_s$ to denote the restriction of $c$ onto $x_s$. Thus $c=\bdd(M,(c_{\{i\}})_{i\le k+1})$ and the sequence $c_{\{0\}},\ldots,c_{\{k+1\}}$ is independent over $M$. Moreover notice that $c_s\equiv_M b_s$ and also $c_s=\bdd(M,(c_{\{i\}})_{i\in s})$. In particular, for $t=\{0,\ldots,k\}$ we have that $c_t\equiv_M b_t$ and so we can find some $d$ such that $c_tc_{\{k+1\}}\equiv_M b_t d$. Hence $d$ is independent from $b_{\{0\}},\ldots,b_{\{k\}}$ over $M$ and in particular, the partial type
$$
\bigcup_{s\in\mathcal P^{-}(k+1)} \tp(d/b_{s})
$$ does not fork over $M$. Furthermore, for any $s\in\mathcal P^{-}(k+1)$ we also have that $c_{s\cup\{k+1\}}\equiv_M b_{s\cup\{k+1\}}$ since the set $s\cup\{k+1\}$ belongs to $\mathcal P^{-}(k+2)$, thus we get $c_sc_{\{k+1\}}\equiv_M b_sb_{\{k+1\}}$ and so
$$
b_{s}d\equiv_M c_{s} c_{\{k+1\}} \equiv_M b_sb_{\{k+1\}}.
$$
Therefore, the partial type
$$
\bigcup_{s\in\mathcal P^{-}(k+1)} \tp(b_{\{k+1\}}/b_{s})
$$
does not fork over $M$, and even less does
$$
\bigcup_{s\in\mathcal P^{-}(k+1)} \tp(a_{k+1}/M,(a_j)_{j\in s}).
$$
Hence, as the language is relational of arity at most $k+1$ and each subset $s$ from $\mathcal P^{-}(k+1)$ has size at most $k$, elimination of quantifiers yields that the latter partial type determines the complete type $\tp(a_{k+1}/M,a_0,\ldots,a_{k})$. Thus, the element $a_{k+1}$ is independent from $a_0,\ldots,a_{k}$ over $M$ and therefore $a_0,\ldots,a_{k+1}$ are independent over $A$ by transitivity, as desired.\qed

By \cite[Proposition 3]{Poi} in the supersimple framework $k$-triviality and triviality agree and therefore:

\begin{cor}\label{CorSupersimpleTrivial}
A supersimple theory with elimination of quantifiers in a relational language of arity at most $k+1$ is trivial whenever it satisfies the $(k+2)$-complete amalgamation property over models.
\end{cor}
Furthermore, as all stable theories have $n$-complete amalgamation over models for every natural number $n$, we immediately obtain that a stable (superstable) theory with elimination of quantifiers in a relational language of arity at most $k+1$ is $k$-trivial (respectively, trivial).

Next we recall a stronger notion of triviality which seems to be more appropriate for stable theories.

\begin{defn}
A theory is {\em $k$-totally trivial} if for every set $A$ and any tuples $a_0,\ldots,a_{k+1}$ with $a_{k+1}$ independent over $A$ from every $k$-tuple formed by the $a_0,\ldots,a_{k}$ we have that $a_{k+1}$ is independent from all of them over $A$.
\end{defn}

It is clear that $k$-totally trivial theories are $k$-trivial and in fact for $k=1$ both notions of triviality agree when the ambient theory has finite $\SU$-rank, see \cite[Proposition 5]{Poi}. 


\begin{prop}\label{PropTotTriv}
A stable theory with elimination of quantifiers in a relational language of arity at most $k+1$ is $k$-totally trivial.
\end{prop}
\pf As pointed out in \cite{Poi}, it is enough to show that its theory is totally trivial among real tuples.
Let $a_0,a_1,\ldots,a_{k+1}$ be  tuples from the home sort and suppose that $a_{k+1}$ is independent from any $k$-tuple from $a_0,\ldots,a_{k}$ over $A$. After replacing $A$ by a model independent from $a_0,\ldots,a_{k+1}$ over $A$, we may assume that $A$ is model. Thus, by stationarity of $\tp(a_{k+1}/A)$ we get that
$$
\bigcup_{i\le k}\tp(a_{k+1}/A,(a_j)_{j\neq i})
$$
does not fork over $A$. Hence, elimination of quantifiers yields that this partial type determines indeed the complete type $\tp(a_{k+1}/A,a_0,\ldots,a_{k})$ since the language is relational of arity at most $k+1$, and therefore $a_{k+1}$ is independent from $a_0,\ldots,a_{k}$ over $A$.\qed

In contrast with the trivial case, there is a superstable theory which is $k$-totally trivial but not totally trivial, see \cite{Poi}. Nevertheless, in the finite $\SU$-rank case all notions agree.

\begin{lemma}
A supersimple $k$-totally trivial theory of finite $\SU$-rank is totally trivial.
\end{lemma}
\pf As it is $k$-totally trivial it is clearly $k$-trivial and so trivial as the ambient theory is supersimple. Moreover, as it has finite $\SU$-rank triviality implies total triviality. \qed

\begin{cor}
A superstable theory of finite $\SU$-rank with elimination of quantifiers in a relational language of bounded arity is totally trivial.
\end{cor}
Next we present some examples of $\omega$-stable totally trivial theories which have quantifier elimination (after expanding the language) in a binary relational language. After Proposition \ref{PropTotTriv}, totally triviality can be seen as a limitation of the language.

\begin{expl}
The free pseudoplane is a bicolored infinite branching graph with no loops. Its theory is $\omega$-stable and admits quantifier elimination after adding for each natural number $n$ the binary relation $d_n(x,y)$ interpreting that ``the distance between $x$ and $y$ is exactly $n$.'' This theory is a well-known example of an infinite Morley rank $\omega$-stable totally trivial theory which is not one-based.
\end{expl}

A more elaborate family of examples which encloses the free pseudoplane is given in \cite{BMPZ}. More precisely, a complete first-order theory is associated to a given a right-angled building with infinite residues.

\begin{expl}
Given a finite graph $\Gamma$, a $\Gamma$-graph is a colored graph with colors $\mathcal A_\gamma$ for $\gamma$ in $\Gamma$, which has no edges between elements of colors $\mathcal A_\gamma$ and $\mathcal A_\delta$ if $\gamma$ and $\delta$ are not adjacent. A flag $F$ of the $\Gamma$-graph is a subgraph $F =\{f_\gamma\}_{\gamma\in \Gamma}$, where each $f_\gamma$ has color $\mathcal A_\gamma$, such that the map $\gamma\mapsto f_\gamma$ induces a graph isomorphism between $\Gamma$ and $F$.
Moreover, given a fixed subset $A$ of $\Gamma$, two flags $F_1$ and $F_2$ are $A$-equivalent, denoted by $F_1\sim_A F_2$, if the set of colors where they differ is contained in $A$. It is clear that this is an equivalence relation. When two flags $F$ and $G$ are $A$-equivalent it is possible to obtain a finite sequence $F_0=F,F_1,\ldots,F_n=G$ of flags such that the colors where $F_i$ and $F_{i+1}$ differ form a non-empty connected subset of $A$. In fact, this path of flags can be obtained in a {\em reduced} way. Due to the technicalities of the definitions we avoid to introduce the precise definition here. We refer the interested reader to \cite[Section 3]{BMPZ} for details.

The theory ${\rm PS}_\Gamma$ in the language of graphs with unary predicates for the colors $\{\mathcal A_\gamma\}_{\gamma\in\Gamma}$ is axiomatized by a collection of sentences expressing that the structure is a $\Gamma$-graph satisfying:
\begin{enumerate}
\item it is a $\Gamma$-space, i.e. every vertex belongs to a flag, and any two adjacent vertices can be expanded to a flag;
\item it is simply connected: there are no non-trivial closed reduced paths between flags;
\item for any $\gamma$ in $\Gamma$, the $\sim_\gamma$-equivalence classes are infinite.
\end{enumerate}
The free pseudoplane corresponds to the theory ${\rm PS}_\Gamma$ with $\Gamma$ being the complete graph $K_2$. The theory ${\rm PS}_\Gamma$ is $\omega$-stable of infinite Morley rank, and it is bi-interpretable with the theory of the induced structure on its space of flags. The latter admits quantifier elimination after adding to the language,  for each reduced word $u$, the definable binary relation $P_u(X,Y)$ defining that ``the flags $X$ and $Y$ are connected by a path of word $u$.'' This is \cite[Theorem 7.24]{BMPZ}. Therefore, total triviality of the theory on the space of flags can be explained by Proposition \ref{PropTotTriv}. Consequently, since total triviality is preserved under interpretation, it follows that the theory ${\rm PS}_\Gamma$  is totally trivial, see \cite[Proposition 7.26]{BMPZ}.
\end{expl}

All these previous examples of superstable theories are totally trivial but not one-based. Nevertheless, as far as $\omega$-categorical structures are concerned, triviality implies one-basedness and finite $\SU$-rank whenever the structure is supersimple.

\begin{lemma}
A finitely based regular type in a countable $\omega$-categorical simple trivial theory is non-orthogonal to a $\SU$-rank one type based over the same parameters.
\end{lemma}
\pf Let $p$ be a regular type (or of pre-weight one) which we may assume to be defined without parameters. By assumption, the relation $E$ defined among realizations of $p$ by
$$
xEy \ \Leftrightarrow \ x\nind y
$$
is an equivalence relation;  moreover, it is $\emptyset$-definable by $\omega$-categoricity. Now, fix some realization $a$ of $p$ and consider its $E$-equivalence class $a_E$. We show that the type $\tp(a_E)$ has $\SU$-rank one. To do so, suppose towards a contradiction that $a_E\nind A$ but $a_E$ is not algebraic over $A$. Note that the representatives of any two distinct realizations of $\lstp(a_E/A)$ are independent over $\emptyset$ as they are not $E$-related and so, by triviality any Morley sequence in $\lstp(a_E/A)$ is also independent over $\emptyset$. Thus $\cb(a_E/A)$ is algebraic  and hence $a_E$ is independent from $A$, a contradiction. Therefore, the result follows since $p=\tp(a)$ is clearly  non-orthogonal to $\tp(a_E)$. \qed

\begin{cor}\label{CorFinite}
A countable $\omega$-categorical supersimple trivial theory has finite $\SU$-rank and so it is one-based.
\end{cor}
\pf Noticing that a finitely based type of infinite monomial $\SU$-rank $\omega$ would be orthogonal to any type of finite $\SU$-rank, we deduce that the ambient theory has finite $\SU$-rank by the previous lemma. Moreover, it is one-based by \cite[Proposition 9]{Poi}. \qed


\section{Finitely homogeneous structures}

In this last section we prove that $2$-transitive finitely homogeneous supersimple structures, satisfying the $n$-complete amalgamation property for large enough $n$, are random structures. Next we give the definition of a random structure, which weakens the definition given in \cite[Definition 2.1]{Koponen3}. Nevertheless,  both notions agree for binary languages or more generally when the language has only relations of a fixed arity.

\begin{defn}\label{DefRandom}
Let $\M$ be a structure in a finite relational language $\mathcal L$. We say that $\M$ is a {\em random structure} if there is no finite $\mathcal L$-structure $A$ for which every proper $\mathcal L$-substructure is embeddable in $\M$ but the structure given as the free amalgamation of all proper $\mathcal L$-substructures of $A$ is not embeddable.
\end{defn}

The definition given by Koponen is stated in terms of minimal forbidden configurations with respect to reducts of $\M$, whereas in the definition above the whole quantifier-free type of every proper substructure is taken into account. In fact, observe that if the size of $A$ is greater than the arity of any relation in the language, the above simple says that $A$ is a minimal forbidden configuration.
Moreover, we point out that one can see randomness as an amalgamation property among quantifier-free types. In particular, it is a property of the theory.

\begin{remark}\label{RemRandom}
A relational structure is random if and only if for any natural number $n$, given a compatible family of  quantifier-free complete types $\{\pi_s((x_i)_{i\in s})\}_{s\in\mathcal P^{-}(n)}$ with $|x_i|=1$ and such that each $\pi_s((x_i)_{i\in s})$ is consistent, we have that $$\bigcup_{s\in\mathcal P^{-}(n)} \pi_s((x_i)_{i\in s})$$ is consistent as well.
\end{remark}


\begin{defn}
A structure $\M$ is said to be {\em primitive} if there is no $\emptyset$-invariant equivalence relation among elements of $\M$. Equivalently, there is a unique $1$-type without parameters which in addition is strong. We say that a structure is {\em $2$-transitive} if there is a unique $2$-type among pairs of distinct elements.
\end{defn}

The following result is due to Aranda L\'opez \cite[Proposition 3.3.3]{Andres}.

\begin{fact}
A primitive binary finitely homogeneous simple structure of $\SU$-rank one is random.
\end{fact}
Koponen was able to extend this result to one-based structures by showing the following in \cite{Koponen3}.
\begin{fact}
A primitive binary finitely homogeneous simple one-based structure has $\SU$-rank one and so it is random.
\end{fact}

Firstly, we shall see that Koponen's result is sharp. This is exhibited in the following example from \cite[Example 3.3.2]{Dugald}, whose existence was pointed out to us by Koponen.

\begin{expl}
Let $\Gamma$ be the graph whose vertices are the $2$-sets of $\mathbb N$ with the relation $S(x,y)$ defined among the $2$-sets if their intersection is a singleton. In addition we add a ternary relation $R(x,y,z)$ to the language interpreted as $x$, $y$ and $z$ are three distinct $2$-sets all sharing a single element. This structure is a primitive finitely homogeneous $\omega$-stable structure whose unique $1$-type has $\SU$-rank two. Moreover, the general theory (or an easy argument) yields that it is one-based and indeed trivial. On the other hand, as the formula $S(x,y)$ $3$-divides over $\emptyset$, the structure $\Gamma$ cannot be a random structure.
\end{expl}

It is worth mentioning that in the example above there are imaginary elements which are not eliminable; namely, the canonical parameters of sets codifying single elements. In fact, this is the unique obstacle to generalize Koponen's result to arbitrary languages.

\begin{lemma}\label{Lemma1}
A countable $\omega$-categorical primitive structure with a simple one-based theory is supersimple of $\SU$-rank one whenever it geometrically eliminates imaginaries.
\end{lemma}
\pf We shall be working inside a large saturated and homogeneous model $\M$ of the given structure.

First, we observe that the algebraic closure of an element of the monster model $\M$ of the theory is a singleton. For this, consider the following equivalence relation defined on $\M\times \M$:
$$xRy \ \Leftrightarrow \ \acl(x)=\acl(y),$$
which is clearly $\emptyset$-invariant, and so $\emptyset$-definable by $\omega$-categoricity. Moreover, as there is
no non-trivial $\emptyset$-invariant  equivalence relation defined on $\M$, the relation $R$ is trivial. Hence, either any class has a single element or there is only one class.
However, the latter would imply that all elements have the same algebraic closure and so $\M$ would be the algebraic closure of a single element, a contradiction. Therefore, every class has a single element. In fact, note that the algebraic closure of a single element is finite by $\omega$-categoricity and thus $|\acl(x)|=|\acl(y)|$ for any two elements since there is only a unique $1$-type without parameters as our structure is primitive. This yields that $\acl(x)=\{x\}$ for any element $x$ of $\M$.

Now, we prove that an arbitrary $1$-type does not fork over the empty-set whenever it is not algebraic. Let $a$ be an element of $\M$ and let $\bar b$ be a tuple of elements of $\M$ such that $a$ is not algebraic over $\bar b$. Consider the canonical base $\cb(\bar b/a)$ and note that it is interalgebraic with a real tuple $\bar c$ by assumption, possibly empty. If $\bar c$ were non-empty, then it would be algebraic over $a$, and so $\bar c= a$. However, this would yield that $a\in \acl(\bar b)$ by one-basedness, a contradiction. Therefore $\cb(\bar b/a)$ belongs to $\acl^{\rm eq}(\emptyset)$ and so $a$ is independent from $\bar b$ by symmetry, as desired. \qed

Now we show our main result, which applies to finitely homogeneous structures without binary relations in the language. Thus, this can be seen as an orthogonal statement to the aforementioned result of  Koponen.

\begin{theorem}\label{ThmRandom}
A $2$-transitive finitely homogeneous supersimple structure in a relational language of arity at most $n$ and whose theory has $(n+1)$-complete amalgamation is a random structure.
\end{theorem}
\pf The theory of a such finitely homogeneous supersimple structure is trivial by Corollary \ref{CorSupersimpleTrivial}. Thus, it is one-based of finite $\SU$-rank by Corollary \ref{CorFinite} and so it is indeed totally trivial.

First, we claim that any two distinct elements are independent: By simplicity there are two distinct elements $a$ and $b$ in the monster model $\M$ which are independent. Thus, as there is a unique $2$-type any two elements must be independent.

Now, as the theory is totally  trivial,  the above claim  yields that any two sets are independent over their intersection and so the algebraic closure of a set is the set itself.  On the other hand, by totally triviality we get for a type $\tp(\bar a/M)$ that
$$
\cb(\bar a/M)=\dcl^{\rm eq}\big( (\cb(a/M))_{a\in \bar a} \big).
$$
Thus, we may distinguish two cases. Suppose first that an element $a\in\bar a$ belongs also to $M$, so  the canonical base $\cb(a/M)=\cb(a/a)$ is definable over $a$ and $a\in\acl(\cb(a/M))$. Thus, as $\tp(a/\cb(a/M))$ is Lascar strong, it implies the type $\tp(a/a)$, hence $a$ is definable over $\cb(a/M)$ and whence $\cb(a/M)=\dcl^{\rm eq}(a)$. On the other hand, if $a$ does not belong to $M$, then $a$ is independent from $M$ and so $\cb(a/M)$ belongs to $\dcl^{\rm eq}(\emptyset)$ since the type $\tp(a)$ is Lascar strong by assumption. Therefore, we obtain that $\cb(\bar a/M)$ is interdefinable with $\bar a\cap M$ and in particular, it is eliminable. Consequently,  by a standard argument, it follows that the theory has weak elimination of imaginaries. Namely, let $e$ be an imaginary given as the equivalence class of a  real tuple $\bar a$ modulo an $\emptyset$-definable equivalence relation, and observe that $e\in\acl^{\rm eq}(\cb(\bar a/e))$. Thus, the type $\tp(\bar a/\cb(\bar a/e))$ implies $\tp(\bar a/e)$ and so for any automorphism $f\in \Aut(\M)$ fixing $\cb(\bar a/e)$ we have that $\bar a\equiv_e f(\bar a)$. Hence $e$ is also the equivalence class of $f(\bar a)$ and so $f$ fixes $e$. Whence $e\in\dcl^{\rm eq}(\cb(\bar a/e))$. On the other hand, the canonical base $\cb(\bar a/e)$ is interdefinable with a real tuple $\bar a'$ and so $e\in\dcl^{\rm eq}(\bar a')$ and $\bar a'\subseteq \acl(e)$, as desired. Therefore, for any real set $A$ we obtain that
\begin{equation}\label{eq}
\acl^{\rm eq}(A)=\dcl^{\rm eq}(\acl(A))=\dcl^{\rm eq}(A). \tag{$\dagger$}
\end{equation}

Finally,  to prove that the given structure is random we use Remark \ref{RemRandom}. To do so, consider a family of compatible quantifier-free complete types $\{\pi_s((x_i)_{i\in s})\}_{s\in\mathcal P^{-}(m)}$ such that each of them is consistent with the ambient theory and each variable $x_i$ has length one. Notice that each quantifier-free type determines a complete type by elimination of quantifiers. Hence, if $m\le n+1$ then, the family $\{\pi_s((x_i)_{i\in s})\}_{s\in\mathcal P^{-}(m)}$ determines an independent system of strong types by (\ref{eq}) and whence, the union of all the $\pi_s$ is consistent by $m$-complete amalgamation.
Otherwise, in case $m>n+1$, fix a realization $\bar c=(c_i)_{i\in t}$ of $\pi_t$ for $t=\{n+1,\ldots,m-1\}$ and note that each quantifier-free type $\pi_s((x_i)_{i\in s\setminus t},(c_i)_{i\in t\cap s})$ is consistent by invariance. Now, for a set $u$ in $\mathcal P^{-}(n+1)$, set $\pi_u'(x_u)$ to be $\pi_{u\cup t}((x_i)_{i\in u},\bar c)$, and observe that the family $\{\pi_u'(x_u)\}_{u\in\mathcal  P^{-}(n+1)}$ determines an independent system of strong types over $\bar c$ by (\ref{eq}) . Therefore, the $(n+1)$-complete amalgamation yields the consistency of the union of all $\pi_u'$ and so  the set
$$
\bigcup_{s\in\mathcal P^{-}(m)} \pi_s((x_i)_{i\in s})
 $$
is consistent. This finishes the proof. \qed

\begin{remark}\label{RemPrimitive}
An inspection of the proof together with Lemma \ref{Lemma1} yields  that any primitive finitely homogeneous supersimple structure in a relational language of arity at most $n$ and whose theory has $(n+1)$-complete amalgamation and eliminates imaginaries geometrically is a random structure. In that case, the theory has $\SU$-rank one by Lemma \ref{Lemma1} and so one can see that any two elements are independent since the algebraic closure of a singleton is the singleton itself. Hence, we obtain that any two sets are independent over their intersection and so the algebraic closure of a set is the set itself. The rest of the proof goes through.
\end{remark}

\end{document}